\documentclass[11pt]{article}
\usepackage{latexsym,amsmath,amsfonts,graphics}
\normalsize
\begin{document}
\begin{center}
{\LARGE\textbf{The Distant-$l$ Chromatic Number \\of Random Geometric Graphs}}\\
\bigskip
Yilun Shang\footnote{Department of Mathematics, Shanghai Jiao Tong University, Shanghai 200240, CHINA. email: \texttt{shyl@sjtu.edu.cn}}\\

\end{center}

\begin{abstract}
A random geometric graph $G_n$ is given by picking $n$ vertices in
$\mathbb{R}^d$ independently under a common bounded probability
distribution, with two vertices adjacent if and only if their
$l^p$-distance is at most $r_n$. We investigate the distant-$l$
chromatic number $\chi_l(G_n)$ of $G_n$ for $l\ge1$. Complete
picture of the ratios of $\chi_l(G_n)$ to the chromatic number
$\chi(G_n)$ are given in the sense of almost sure convergence.
\bigskip

\textbf{Keywords:} chromatic number; distance coloring; random
geometric graph.
\end{abstract}

\bigskip
\normalsize

\noindent{\Large\textbf{1. Introduction}}

\smallskip
We consider in this short paper the distance coloring of random
geometric graph $G_n:=G(\mathcal{X}_n,r_n)$, which is obtained as
follows. Take $l^p$-norm $||\cdot||$ in $\mathbb{R}^d$ for $1\le
p\le\infty$. Let $f$ be some bounded probability density function on
$\mathbb{R}^d$ and let $\mathcal{X}_n=\{X_1,X_2,\cdots,X_n\}$, where
$\{X_i\}$ are i.i.d. random $d$-vectors with the common density $f$.
Let $r_n$ be a sequence of distance satisfying $r_n\rightarrow0$ as
$n\rightarrow\infty$. Then we denote by $G(\mathcal{X}_n,r_n)$ the
graph with vertex set $\mathcal{X}_n$ and with edges $X_iX_j$, if
$||X_i-X_j||<r_n$ for all $i\not=j$. An excellent introduction to
random geometric graphs is available in \cite{1}. The vertex
coloring in geometric graphs is closely related with the
\textit{radio channel assignment problem}, see e.g. \cite{16,5,17}
for its history as well as an extensive treatment of this
significant issue.

Recall that a $k$-coloring of a graph $G$ is a map $g:
V(G)\rightarrow \{1,2,\cdots,k\}$ such that $g(u)\not=g(v)$ whenever
$uv\in E(G)$ and that the chromatic number $\chi(G)$ is the least
$k$ for which $G$ is $k$-colorable. For a graph $G$, the graph
distance $d_G(u,v)$ between two vertices $u$ and $v$ is defined as
the length of a shortest path joining them (hence the graph distance
will be infinity if they are in different components). Thereby we
get two kinds of distance between two vertices in our geometric
setting, i.e. $d_G(u,v)$ and $||u-v||$, whose interrelationship
(Lemma 4) turns to be important in the proof of our main theorems.
For $l\ge 1$, a distant-$l$ coloring of $G$ is a coloring of the
vertices such that vertices at distance $(d_G)$ less than or equal
to $l$ have different colors. The least number for which a
distant-$l$ coloring exists is  called the \textit{distant-$l$
chromatic number} of $G$, designated by $\chi_l(G)$. Recall that a
distant-$l$ coloring of $G$ is equivalent to an ordinary vertex
coloring of $l$ power $G^l$ of $G$. (The $l$ power of a graph $G$,
denoted as $G^l$, is the graph with the same vertex set and in which
two vertices are joined by an edge if and only if they have distance
$(d_G)$ less than or equal to $l$ in $G$.) Hence
$\chi(G^l)=\chi_l(G)$, particularly, $\chi_1(G)=\chi(G)$.

Distance coloring has been a long standing topic in graph theory and
has been dealt with mostly in planar graphs. We refer the reader to
\cite{9,7,8} and \cite{18} for more details regarding this subject.
Recently, D\'{\i}az et al. \cite{6} studied the distant-2 chromatic
number of $G_n$ in the Euclidean plane by takeing $f=1_{[0,1]^2}$.
Their results show that the order of $\chi_2(G_n)$ is consistent
with that of $\chi(G_n)$ in the connectivity regime, i.e. when
$nr_n^2=\Theta(\ln n)$. In the study of random geometric graph
$G_n$, some limiting regimes for $r_n$ are of special interest
\cite{1}. One of these is connectivity regime in which
$nr_n^d=\Theta(\ln n)$ as mentioned above. When $nr_n^d\gg\ln n$ we
refer to the limiting regime as the superconnectivity regime and the
limiting regime $nr_n^d\ll\ln n$ is called the subconnectivity
regime. Our aim, in this paper, is to determine the strong law
results of the ratios of $\chi_l(G_n)$ to $\chi(G_n)$ in all the
above three cases. In addition, a focusing phenomenon that the
probability measure becomes concentrated on two consecutive integers
is observed in \cite{2} for $\chi(G_n)$ in the subconnectivity
regime. We shall state an analogous result (Theorem 1) for
$\chi_l(G_n)$, whose proof may be adapted from \cite{2}
straightforwardly.

Before going further, we introduce some preliminary definitions. Let
$\mathrm{vol}(\cdot)$ denote the $d$-dimensional Lebesgue measure.
In the rest of the paper, let
$f_{\max}:=\sup\{t|\mathrm{vol}(\{f(x)>t\})>0\}$ be the essential
supremum of the probability density function $f$. For any graph $G$,
we denote the maximum degree of $G$ by $\triangle(G)$ and the clique
number of $G$ by $\omega(G)$. Recall that we have the basic
inequalities: $\omega(G)\le\chi(G)\le\triangle(G)+1$. Given
$x\in\mathbb{R}^d$ and $r>0$, let $B(x,r)$ be the ball centered at
$x$ with radius $r$. We will say that a sequence of events $E_n$
holds a.s. (almost surely) if $P(E_n)\rightarrow1$ as
$n\rightarrow\infty$.

The rest of this paper is organized as follows. In Section 2, we
give our main results for the distant chromatic number. Section 3
contains the proofs. We conclude the paper in Section 4.

\bigskip
\noindent{\Large\textbf{2. Statement of main results}}
\smallskip

Throughout the paper, we assume the probability density function $f$
is bounded, that is, $f_{\max}<\infty$. We alluded to the following
result in Section 1.

\medskip
\noindent\textbf{Theorem 1.(Focusing)} \itshape If $nr_n^d=o(\ln n)$
and $l\ge1$ fixed, then there exists a sequence $\{a_n\}$ such that
$$P(a_n\le\chi_l(G_n)\le a_n+1)\rightarrow1,\qquad as\ \ n\rightarrow\infty.$$
\normalfont

We leave the proof as an exercise for the readers consulting
Corollary 2 of \cite{2}.

Let $\mathrm{supp}f$ be the support of $f$, i.e.
$\mathrm{supp}f:=\overline{\{x|f(x)>0\}}$. Let $f_0$ be the
essential infimum of $f$ over $\mathrm{supp}f$, that is, the largest
$h$ such that $P(f(X_1)\ge h)=1$.

\medskip
\noindent\textbf{Theorem 2.(Superconnectivity regime)} \itshape
Suppose $nr_n^d\gg\ln n$ and $l\ge1$ fixed. Let $f$ satisfies (a)
a.e. continuous and $\mathrm{supp}f$ does not contain a sequence of
isolated points which has a limit point in $\mathrm{supp}f$; or (b)
$f_0>0$, then
$$
\frac{\chi_l(G_n)}{l^d\chi(G_n)}\rightarrow1\qquad a.s.
$$
\normalfont

We remark that the conditions (a) and (b) imposed on $f$ above are
rather mild; in fact, typical distributions such as normal
distribution and $f=1_{[0,1]^d}$ are clearly allowed.

\medskip
\noindent\textbf{Theorem 3.(Connectivity regime)} \itshape Suppose
$nr_n^d=\Theta(\ln n)$ and $l\ge1$ fixed, then
$$
\frac{\chi_l(G_n)}{l^d\chi(G_n)}\rightarrow c\qquad a.s.
$$
where $c\in [l^{-d},1]$ is a constant depending only on the quantity
involved in ``$\Theta$''.\normalfont\medskip

The constant $c$ will be explicitly given in the proof.

\medskip
\noindent\textbf{Theorem 4.(Subconnectivity regime)} \itshape
Suppose $n^{-\varepsilon}\ll nr_n^d\ll \ln n$ for all
$\varepsilon>0$ and $l\ge1$ fixed, then
$$
\frac{\chi_l(G_n)}{\chi(G_n)}\rightarrow1\qquad a.s.
$$
\normalfont

Notice that Theorem 2 and 3 rely explicitly on the dimension of the
underlying space $\mathbb{R}^d$ while Theorem 4 does not.

To close up the spectrum of limiting regimes, we observe (by
exploiting a result in \cite{1} Section 6.1) that if $nr_n^d\le
n^{-\varepsilon}$ for some $\varepsilon>0$, then there exists some
$c>0$ such that $P(\chi_l(G_n)\le c)\rightarrow1$ and
$P(\chi(G_n)\le c)\rightarrow1$ as $n\rightarrow\infty$. Hence there
won't be any interesting strong law in this case.

We refer the readers to \cite{3} for a number of results regarding
the strong law of large numbers in chromatic number $\chi(G_n)$,
which are largely improved than those discovered earlier by Penrose
et al. (see e.g.\cite{1}). Wherefore our theorems suggest the strong
laws of $\chi_l(G_n)$.

\bigskip
\noindent{\Large\textbf{3. Proofs}}

\smallskip

We will need some strong law results from \cite{3}, which take an
important role in the proofs. To make the present work
self-contained, some technical definitions are included as follows.
For a measurable set $A\subseteq\mathbb{R}^d$, if
$\lim_{\varepsilon\rightarrow0}\mathrm{vol}(A_{\varepsilon})=\mathrm{vol}(A)$,
where $A_{\varepsilon}:=A+B(0,\varepsilon)$, then we say $A$ has a
small neighborhood. Let $\mathcal{F}$ be the collection of all
non-negative, bounded, measurable functions
$\varphi:\mathbb{R}^d\rightarrow\mathbb{R}$ with
$0<\mathrm{vol}(\mathrm{supp}\varphi)<\infty$, and
$\{x|\varphi(x)>a\}$ having a small neighborhood for all
$a\in\mathbb{R}$. Given $\varphi\in\mathcal{F}$, let
$M_{\varphi}:=\sup_{x\in\mathbb{R}^d}\sum_{i=1}^{n}\varphi(\frac{X_i-x}{r_n})$.
Set a function $H(x):=x\ln x-x+1$ for $x>0$. For
$\varphi\in\mathcal{F}$ and $0<t<\infty$, let us set
$$\xi(\varphi,t):=\int_{\mathbb{R}^d}\varphi(x)e^{s\varphi(x)}\mathrm{d}x,$$
where $s=s(\varphi,t)$ is the unique non-negative solution to the
equation
$\int_{\mathbb{R}^d}H(e^{s\varphi(x)})\mathrm{d}x=\frac{1}{tf_{\max}}$.
We also set
$\xi(\varphi,\infty)=\int_{\mathbb{R}^d}\varphi(x)\mathrm{d}x$
naturally. Let $\mathcal{G}$ be the collection of measurable,
non-negative functions $\varphi:\mathbb{R}^d\rightarrow[0,1]$ such
that $\sum_{x\in A}\varphi(x)\le 1$ for any set
$A\subseteq\mathbb{R}^d$ satisfying $||x-y||>1$ for all $x\not=y\in
A$. Denote $k_n:=\frac{\ln n}{\ln(\ln n/nr_n^d)}$ throughout the
paper.

\medskip
\noindent\textbf{Lemma 1} \itshape(\cite{3}). Suppose
$\frac{nr_n^d}{\ \ln n\ }\rightarrow t\in(0,\infty]$, then
$$
\frac{\chi(G_n)}{\ nr_n^d\ }\rightarrow
f_{\max}\sup_{\varphi\in\mathcal{G}}\xi(\varphi,t) \qquad a.s.
$$
\normalfont\smallskip

\noindent\textbf{Lemma 2} \itshape(\cite{3}). Let
$W\subseteq\mathbb{R}^d$ be a bounded, measurable set with non-empty
interior and having a small neighborhood. Suppose $\varphi=1_W$. For
every $\varepsilon>0$, there exists a $\delta>0$ such that if
$n^{-\delta}\le nr_n^d\le\delta\ln n$ then
$$
P((1-\varepsilon)k_n\le M_{\varphi}\le(1+\varepsilon)k_n\ for\ all\
but\ finitely\ many\ n)=1,
$$
where $k_n$ is defined as above.\normalfont\medskip

We will also need the following property for the functional
$\xi(\varphi,t)$.

\medskip
\noindent\textbf{Lemma 3} \itshape(\cite{3}). For $t,h>0$ and
non-negative, bounded, measurable, integrable function $\varphi$, we
have
$$\Big(\frac{t}{t+h}\Big)\xi(\varphi,t)\le\xi(\varphi,t+h)\le\xi(\varphi,t).$$
\normalfont\smallskip

The next lemma reveals that the shortest path between any pair of
nodes in $G_n$ is close to a straight line in and above the
connectivity regime. This result can be improved, but it will be
enough for our purpose here.\bigskip

\vspace{10mm} \noindent\textbf{Lemma 4.} \itshape Suppose
$nr_n^d\ge\Theta(\ln n)$ and $l\ge 1$ fixed. Let $P_l^n$ denote the
probability $P(d_{G_n}(X_1,X_2)>l,||X_1-X_2||<lr_n)$, then
$$
P_l^n=o(n^{-2})\qquad as\ n\rightarrow\infty.
$$
\normalfont\smallskip \noindent\textbf{Proof.} We use the inductive
method for $l$. If $l=1$, then
$$
P_1^n=P(d_{G_n}(X_1,X_2)>1,||X_1-X_2||<r_n)=0, \qquad\forall\ n\ge
2.
$$
Now assuming that $P_l^n=o(n^{-2})$ for some $l\ge1$, we aim to
prove $P_{l+1}^n=o(n^{-2})$. Assume that
$||X_1-X_2||=(l+1-\varepsilon)r_n$ with $0<\varepsilon\le1$. Let
$Q:=B(X_1,(l-\varepsilon/2)r_n)\cap B(X_2,r_n)$, see e.g. Fig.1
(which is drawn under $l^2$-norm and $d=2$). We have
\begin{eqnarray}
1-P_{l+1}^n&\ge&1-P(d_{G_n}(X_1,X_2)>l+1\big|\ ||X_1-X_2||<(l+1)r_n)\nonumber\\
&=&P(d_{G_n}(X_1,X_2)\le l+1\big|\ ||X_1-X_2||<(l+1)r_n)\nonumber\\
&\ge&P(\exists X_i\in Q,d_{G_n}(X_1,X_i)\le
l\big|\ ||X_1-X_2||<(l+1)r_n)\nonumber\\
&=&P(\exists X_i\in Q\big|\ ||X_1-X_2||<(l+1)r_n)\nonumber\\
 & &\cdot P(d_{G_n}(X_1,X_i)\le
l\big|\ \exists X_i\in Q, ||X_1-X_2||<(l+1)r_n)\nonumber\\
&:=&P_1(n)\cdot P_2(n)\label{1}
\end{eqnarray}
Fig.1 implies that the value of $r_n$ is only a scaling factor and
that $\mathrm{vol}(Q)$ is proportional to $r_n^d$ (w.r.t. any
$l^p$-norm $||\cdot||$). Moreover, $\mathrm{vol}(Q)=c_1r_n^d$ for
some positive constant $c_1=c_1(l,\varepsilon)$ depending only on
$l$ and $\varepsilon$. Therefore, by the requirement of probability
density $f$ and the asymptotic behavior of $r_n$, we get
$$
\lim_{n\rightarrow\infty}(1-P_1(n))=\lim_{n\rightarrow\infty}\Big(1-\int_Qf(x)\mathrm{d}x\Big)^n\le\lim_{n\rightarrow\infty}e^{-\mathrm{vol}(Q)c_2}=0,
$$
for some positive constant $c_2$. Thus $P_1(n)\rightarrow1$ as
$n\rightarrow\infty$.

On the other hand, by the inductive assumption, $1-P_2(n)\le
P_l^n=o(n^{-2})$. Hence, $P_2(n)\ge 1-o(n^{-2})$ as
$n\rightarrow\infty$. Taking limit in both side of (\ref{1}) gives
$$
P_{l+1}^n=o(n^{-2})\qquad \mathrm{as}\ n\rightarrow\infty,
$$
which concludes the proof. $\Box$ \bigskip

Recall that $G_n=G(\mathcal{X}_n,r_n)$ and let
$G'_n=G(\mathcal{X}_n,lr_n)$, then it's easy to see that
$\chi(G'_n)\ge\chi_l(G_n)\ge\chi(G_n)$. Now for $t\in(0,\infty]$,
suppose $\frac{nr_n^d}{\ \ln n\ }\rightarrow t$ as
$n\rightarrow\infty$. By applying Lemma 1 to graph $G'_n$, and using
Lemma 3, we have
\begin{equation}
\frac{\chi(G'_n)}{l^dnr_n^d}\sim
f_{\max}\sup_{\varphi\in\mathcal{G}}\xi(\varphi,l^dt)\le
f_{\max}\sup_{\varphi\in\mathcal{G}}\xi(\varphi,t)\sim\frac{\chi(G_n)}{nr_n^d}\quad
a.s.\label{2}
\end{equation}
Hence $\chi(G'_n)\le l^d\chi(G_n)$ almost surely for large enough
$n$. Therefore we get
$$
P(\chi(G_n)\le\chi_l(G_n)\le l^d\chi(G_n)\ \mathrm{for}\
\mathrm{all}\ \mathrm{but}\ \mathrm{finitely}\ \mathrm{many}\ n)=1.
$$

Thus it can be seen from Theorem 2 and 4 that the upper bound for
$\chi_l(G_n)$ is asymptotically attained in the superconnectivity
regime while the lower bound is achieved in the subconnective case.

\medskip
\noindent\textbf{Proof of Theorem 2.} Take $t=\infty$ in (\ref{2}).
From the above discussion and definition of the functional
$\xi(\varphi,t)$, we have $\chi(G'_n)=(1+o(1))l^d\chi(G_n)$ almost
surely for large enough $n$. Then it suffices to prove
$P(\chi(G'_n)>\chi_l(G_n))=o(1)$ as $n\rightarrow\infty$.

Now we have by Lemma 4,
\begin{eqnarray*}
P(\chi(G'_n)>\chi_l(G_n))&\le&P(\ \exists X_i,X_j\ s.t.\
d_{G_n}(X_i,X_j)>l\ \mathrm{and}\ ||X_i-X_j||<lr_n)\\
&\le&n^2P_l^n=o(1).
\end{eqnarray*}
as $n\rightarrow\infty$. The proof is then completed. $\Box$

\medskip
\noindent\textbf{Proof of Theorem 3.} Suppose  $\frac{nr_n^d}{\ \ln
n\ }\rightarrow t\in(0,\infty)$ as $n\rightarrow\infty$. Employing
Lemma 4 in the same way as the above proof suggests that
$\chi_l(G_n)=(1+o(1))\chi(G'_n)$. Therefore by the expression
(\ref{2}),
$$
\frac{\chi_l(G_n)}{l^d\chi(G_n)}\sim\frac{\chi(G'_n)}{l^d\chi(G_n)}
\rightarrow\frac{\sup_{\varphi\in\mathcal{G}}\xi(\varphi,l^dt)}{\sup_{\varphi\in\mathcal{G}}\xi(\varphi,t)}:=c(t),
\qquad a.s.
$$
where, by Lemma 3, $c(t)\in[l^{-d},1]$ as claimed. $\Box$

\medskip
\noindent\textbf{Proof of Theorem 4.} Observe that
$$
\omega(G_n)\le\chi(G_n)\le\chi_l(G_n)\le\triangle(G_n^l)+1\le\triangle(G'_n)+1.
$$
By Lemma 2 and the remarks in \cite{3} (see also \cite{17}), we get
almost surely
$$
\omega(G_n)\sim k_n,\quad \chi(G_n)\sim k_n,\quad
\triangle(G'_n)\sim\frac{\ln n}{\ln (\ln n/l^dnr_n^d)}\sim k_n
$$
as $n\rightarrow\infty$. Recall that $k_n=\frac{\ln n}{\ln(\ln
n/nr_n^d)}$ . Also note that $k_n$ tends to infinity by the assumed
asymptotic behavior of $r_n$, which concludes the proof. $\Box$

\bigskip
\noindent{\Large\textbf{4. Concluding remarks}}
\smallskip

We have investigated in this paper the asymptotic behavior of
$\chi_l(G_n)$ when the parameter $l$ is fixed. It is, however,
possible to generalize the results to growing $l$ as long as $l$
dose not increase too quickly. Another issue which we have not
studied but might be of significance in practice is the rates of
convergence of the ratio $\chi_l(G_n)/\chi(G_n)$.

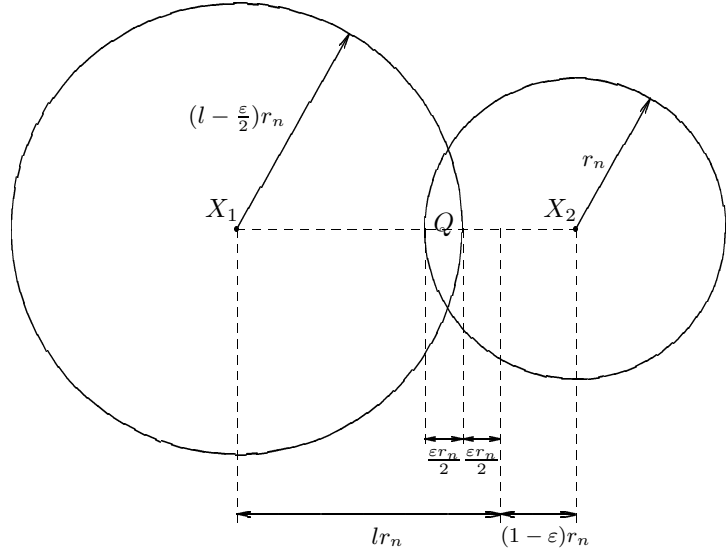
\end{document}